\documentclass[a4paper,12pt]{amsart}
\usepackage{palatino}

\newtheorem{theo}{Theorem}[section]

\newtheorem{fact}[theo]{Fact}

\newtheorem{prop}[theo]{Proposition}
\newtheorem{lem}[theo]{Lemma}

\newtheorem{conj}[theo]{Conjecture}
\newtheorem{defi}[theo]{Definition}

\newtheorem{rema}[theo]{Remark}
\newtheorem{remas}[theo]{Remarks}

\def \kbar {{\overline k}}
\def \Xbar {{\bar X}}

\def \Romannumeral #1 {\expandafter\uppercase\expandafter {\romannumeral #1} }

\def \P {{\bf P}}
\def \tors {{\rm tors}}

\def \Pic {{\rm {Pic}}}

\def \Gal {{\rm{Gal}}}

\def \Frob{{\rm{Frob}}}
\def \calo {{\mathcal O}}
\def \Spec {{\rm{Spec\,}}}

\def \End {{\rm {End}}}

\def\ov{\overline}

\def \Z {{\bf Z}}
\def \Q {{\bf Q}}
\def \F {{\bf F}}

\def \im {{\rm {Im\,}}}
\def \G {{\bf G}_m}

\def \Het {H_{\mbox{\scriptsize\'et}}}

\def\smallsquare{\vbox{\hrule\hbox{\vrule height 1 ex\kern 1 ex\vrule}\hrule}}
\def\enddem{\hfill \smallsquare\vskip 3mm}

\usepackage{amscd}
\usepackage{amssymb}
\usepackage{amsmath}
\usepackage{pb-diagram}

\everymath{\displaystyle}

\def \tors {{\rm tors}}
\def \id {{\rm id}}

\def \CX{{\mathcal X}}
\def \CY{{\mathcal Y}}
\def \lim {{\rm{lim\,}}}
\def \limit{\mathop{\lim}}
\def \limproj{\limit\limits_{\leftarrow}}

\DeclareFontFamily{U}{wncy}{}
\DeclareFontShape{U}{wncy}{m}{n}{%
   <5>wncyr5%
   <6>wncyr6%
   <7>wncyr7%
   <8>wncyr8%
   <9>wncyr9%
   <10>wncyr10%
   <11>wncyr10%
   <12>wncyr6%
   <14>wncyr7%
   <17>wncyr8%
   <20>wncyr10%
   <25>wncyr10}{}
\DeclareMathAlphabet{\cyrille}{U}{wncy}{m}{n}

\def \R{{\bf R}}

\title{Cohomology and torsion cycles over the maximal cyclotomic extension}
\author{by Damian R\"ossler and Tam\'as Szamuely}
\address{Mathematical Institute, University of Oxford, Andrew Wiles Building, Radcliffe Observatory Quarter, Woodstock Road, Oxford OX2 6GG, United Kingdom}
\email{damian.rossler@maths.ox.ac.uk}
\address{Alfr\'ed R\'enyi Institute of Mathematics, Hungarian Academy of Sciences, Re\'altanoda utca 13--15, H-1053 Budapest, Hungary and Central European University, N\'ador utca 9, H-1051 Budapest, Hungary}
\email{szamuely.tamas@renyi.mta.hu}

\date{\today}
\begin{document}
\maketitle

\begin{abstract}
A classical theorem by K. Ribet asserts that an abelian variety defined over the maximal cyclotomic extension $K$ of a number field has only finitely many torsion points. We show that this statement can be viewed as a particular case of a much more general one, namely that the absolute Galois group of $K$ acts with finitely many fixed points on the \'etale cohomology with $\Q/\Z$-coefficients of a smooth proper $\overline K$-variety defined over $K$. We also present a conjectural generalization of Ribet's theorem to torsion cycles of higher codimension. We offer supporting evidence for the conjecture in codimension 2, as well as an analogue in positive characteristic.
\end{abstract}

\section{Introduction}

The Mordell--Weil theorem asserts that the group of points of an abelian variety over a number field is finitely generated. Since Mazur's pioneering paper \cite{mazur} there have been speculations about extensions of this statement to certain infinite algebraic extensions of $\Q$. For the  cyclotomic $\Z_p$-extension of a number field the question is still open in general. However, the Mordell--Weil rank of an abelian variety can be infinite over the maximal cyclotomic extension of a number field obtained by adjoining {\em all} complex roots of unity (see e.g. \cite{rw}). Therefore the following, by now classical, theorem of Ribet \cite{ribet} is all the more remarkable.

\begin{theo}[Ribet]\label{ribet} Let $k$ be a number field, and $K$ the field obtained by adjoining all roots of unity to $k$ in a fixed algebraic closure $\kbar$. If $A$ is an abelian variety defined over $k$, the torsion subgroup of $A(K)$ is finite.
\end{theo}

Note that Ribet's theorem is specific to the maximal cyclotomic extension of $k$. For instance, finiteness of the torsion subgroup may fail if one works over the maximal abelian extension. Indeed, for an abelian variety of CM type all  torsion points are defined over the maximal abelian extension of the number field. Conversely, Zarhin \cite{zarhin} proved that for non-CM simple abelian varieties finiteness of the torsion subgroup holds over the maximal abelian extension as well.\smallskip

 In the present paper we offer two kinds of generalizations of Ribet's theorem. The first one is cohomological, the second is motivic (and largely conjectural). Let us start with the cohomological statement.

\begin{theo}\label{cohribet} Let $k$ and $K$ be as in Theorem \ref{ribet}, and set $G:=\Gal(\kbar|K)$. Consider a smooth proper geometrically connected variety $X$ defined over $k$, and denote by $\Xbar$ its base change to the algebraic closure $\kbar$.

For all {\bf odd} $i$ and all $j$ the groups $\Het^i(\Xbar, \Q/\Z(j))^G$ are finite.
\end{theo}

Here, as usual, we denote by $\Het^i(\Xbar, \Q/\Z(j))^G$ the subgroup of $G$-invariants in the \'etale cohomology group $\Het^i(\Xbar, \Q/\Z(j))$. The theorem is {\em not} true for even degree cohomology (see Remark \ref{remeven} below). In fact, the twist $j$ does not really play a role in the statement since $G$ fixes all roots of unity.

Note that Theorem \ref{ribet} is the special case $i=j=1$ of the above statement. Indeed, if we apply the theorem to the dual abelian variety $A^*$ of an abelian variety $A$ defined over $k$, the Kummer sequence in \'etale cohomology induces a Galois-equivariant isomorphism between
$\Het^1(\overline A^*, \Q/\Z(1))$ and the torsion subgroup of $\Het^1(\overline A^*, \G)$. But since the N\'eron--Severi group of an abelian variety is torsion free, we may identify this torsion subgroup with that of $\Pic^0 (A^*)(\overline K)=A(\overline K)$.

One may also view Ribet's theorem as a finiteness result about the codimension 1 Chow group of smooth projective varieties. In this spirit we propose the following conjectural generalization.

\begin{conj}\label{conj1} Let $k$ and $K$ be as in Theorem \ref{ribet}, and let $X$ be a smooth proper geometrically connected variety defined over $k$. Denote by $X_K$ the base change of $X$ to $K$.

For all $i>0$ the codimension $i$ Chow groups $CH^i(X_K)$ have finite torsion subgroup.
\end{conj}

Recall that for $X$ itself the Chow groups $CH^i(X)$ are conjecturally finitely generated (a consequence of the generalized Bass conjecture on the finite generation of motivic cohomology groups of regular schemes of finite type over $\Z$ -- see e.g. \cite{kahnhb}, 4.7.1).

As is customary with conjectures concerning algebraic cycles, evidence is for the moment scarce in codimension $>1$. However, we have the following positive result.

\begin{theo}\label{main1} In the situation of Conjecture \ref{conj1}, assume moreover that the coherent cohomology group $H^2_{\rm Zar}(X,\calo_X)$ vanishes. Then the torsion subgroup  of $CH^2(X_K)$  has finite exponent. It is finite if furthermore the $\ell$-adic cohomology groups {\rm  $\Het^3(\Xbar, \Z_\ell)$} are torsion free for all $\ell$.
\end{theo}

All geometric assumptions of the theorem are satisfied, for instance, by smooth complete intersections of dimension $>2$ in projective space. Recall also that for a smooth proper variety defined over a finite number field and satisfying $H^2_{\rm Zar}(X,\calo_X)=0$ the torsion part of $CH^2(X)$ is known to be finite by a theorem first proven by Colliot-Th\'el\`ene and Raskind \cite{ctra2}. We shall adapt their methods to our situation.

Motivated by Mazur's program, one might also speculate about finite generation of Chow groups over the cyclotomic $\Z_p$-extension of a number field. We do not pursue this line of inquiry here.\medskip

We now turn to positive characteristic analogues of Ribet's theorem. In the case where $k$ is the function field of a curve $C$ defined over a finite field $\F$, the role of the maximal cyclotomic extension is played by the base extension $k\overline\F$. As the analogue of Theorem \ref{ribet} is plainly false for an abelian variety that is already defined over $\overline\F$,  we have to impose a non-isotriviality condition.

\begin{theo}[Lang--N\'eron]\label{posribet} Let $k$ be the function field of a smooth proper geometrically connected curve $C$ defined over a finite field $\F$, and set $K:=k\overline\F=\overline\F(C)$. If $A$ is an abelian variety defined over $k$ whose base change $A_K$ has trivial $K|\overline\F$-trace , the torsion subgroup of $A(K)$ is finite.
\end{theo}

This is a consequence of the Lang--N\'eron theorem (\cite{lang}, Chapter 6, Theorem 2; see also this reference for the definition of $K|\overline\F$-trace). In fact, under the assumptions of the theorem the group $A(K)$ is even finitely generated.

When seeking an analogue for higher degree cohomology, we have to find a replacement for the non-isotriviality assumption incarnated in the vanishing of the $K|\overline\F$-trace. We propose the following condition.

\begin{defi}\label{large}\rm Let $k=\F(C)$ and $K=\overline\F(C)$ be as in Theorem \ref{posribet}, and let $\ell$ be a prime different from $p={\rm char}(\F)$. Consider a smooth proper geometrically connected variety $X$ defined over $k$, with base change $\Xbar$ to the separable closure $\kbar$. We say that the cohomology group ${\Het^i(\Xbar, \Q_\ell(j))}$ has {\em large variation} if after finite extension of the base field $\F$ there exists a proper flat morphism $\CX\to C$ of finite type with generic fibre $X$ and two $\F$-rational points $c_1,c_2\in C$ such that the fibres $\CX_{c_1}$, $\CX_{c_2}$ are smooth and the associated Frobenius elements $\Frob_{c_1}, \Frob_{c_2}$  act on $\Het^{i}(\Xbar, \Q_\ell(j))$ with coprime
characteristic polynomials.
\end{defi}

Here, as usual, the action of the $\Frob_{c_r}$ $(r=1,2)$ is to be understood as follows. We pick a decomposition group $D_r\subset\Gal(\kbar|k)$ attached to $c_r$; it is defined only up to conjugacy but this does not affect the definition. The smoothness condition on $\CX_{c_r}$ implies that the inertia subgroup $I_r\subset D_r$ acts trivially on cohomology, hence we have an action of $D_r/I_r=\langle\Frob_{c_r}\rangle$ on $\Het^{i}(\Xbar, \Q_\ell(j))$. Furthermore, it is a consequence of the Weil conjectures that the characteristic polynomials in the above definition have rational coefficients and are independent of the prime $\ell$; hence the definition does not depend on $\ell$. For a relation with triviality of the $K|\F$-trace for abelian varieties, see Remarks \ref{abvarremas}.

Based on the above definition, we offer the following higher degree analogue of Theorem \ref{posribet}.

\begin{prop}\label{cohposribet} Let $k$, $K$  and $X$ be as in the previous definition. Assume moreover that $i>0$ and $j\in \Z$ are such that the cohomology group $\Het^i(\Xbar, \Q_\ell(j))$ has large variation.

Then the group $\Het^i(\Xbar, (\Q/\Z)'(j))^G$ is finite, where $G=\Gal(\kbar|K)$ and $$(\Q/\Z)'(j)=\bigoplus_{\ell\neq p}\Q_{\ell}/\Z_{\ell}(j).$$
\end{prop}

In fact, large variation can only occur for odd degree cohomology (Remark \ref{largevarremas}), so the situation indeed resembles that in Theorem \ref{cohribet}.

Building upon this finiteness statement, we propose the following positive characteristic analogue of Conjecture \ref{conj1}.

\begin{conj}\label{conj2} Let $k$, $K$ and $X$ be as in Definition \ref{large}.
Given $i>0$, assume that  the cohomology group $\Het^{2i-1}(\Xbar, \Q_\ell(i))$ has large variation for $\ell\neq p$, where $p={\rm char}(k)$.

Then the prime-to-$p$ torsion subgroup of $CH^i(X_K)$  is finite.
\end{conj}

We offer the following evidence for the conjecture:

\begin{theo}\label{main2} In the situation of Conjecture \ref{conj2}, assume moreover that $X$ is a projective variety which is liftable to characteristic 0 and for which the coherent cohomology group $H^2_{\rm Zar}(X,\calo_X)$ vanishes.  Under the large variation assumption for $i=2$, the prime-to-$p$ torsion subgroup  of $CH^2(X_K)$ is finite.
\end{theo}

The liftability assumption means that $X$ arises as the special fibre of a smooth projective flat morphism $\CX\to S$, where $S$ is the spectrum of a complete discrete valuation ring of characteristic 0 and residue field $k$. It holds for smooth complete intersections or for varieties satisfying the condition $H^2_{\rm Zar}(X, {\mathcal T}_{X/k})=0$ in addition to $H^2_{\rm Zar}(X,\calo_X)=0$, where ${\mathcal T}_{X/k}$ denotes the tangent sheaf (see \cite{fga}, \S 6 or \cite{illusie}, \S 8.5).\smallskip

Many thanks to Anna Cadoret, Jean-Louis Colliot-Th\'el\`ene, Wayne Raskind and Douglas Ulmer. We are also very grateful to the referee for his insightful comments which in particular enabled us to remove an unnecessary restriction in the statement of Theorem \ref{main2}. The second author was partially supported by NKFI grant No. K112735.

\section{Preliminaries in \'etale cohomology}

This section is devoted to an auxiliary statement, presumably well known to some, which will enable us to reduce the case of infinite torsion coefficients to those of $p$-adic and mod $p$ coefficients.

\begin{prop}\label{ABlemma}
Let $X$ be a smooth proper variety over a field $F$ of characteristic 0. Denote by $\Xbar$ its base change to the algebraic closure $\overline F$, and set $G:=\Gal(\overline F|F)$.
Fix a pair $(i,j)$ of nonnegative integers, and assume the following two conditions hold.

\begin{enumerate}
\item[\rm (A)] $\Het^i(\bar X, {\Q}_p(j))^G=0$ for all primes $p$.
\item[\rm (B)] $\Het^i(\bar X, {\Z}/p{\Z}(j))^G=0$ for all but finitely many primes $p$.
\end{enumerate}
Then the group
$\Het^i(\bar X, {\Q}/{\Z}(j))^G
$ is finite.
\end{prop}

\begin{rema}\label{remABlemma}\rm The proposition also holds over fields of  positive characteristic if we restrict to primes different from the characteristic and replace ${\Q}/{\Z}(j)$ by the direct sum of all ${\Q}_p/{\Z}_p(j)$ for primes different from the characteristic. The reader is invited to make the straightforward modifications in the proofs to follow.

\end{rema}

In the proof of the proposition and in later arguments we shall use the following basic properties of \'etale cohomology.

\begin{fact}\label{gabber}\rm
The groups $\Het^i(\Xbar, \Z/p^r\Z(j))$ are finite for all $r$. Moreover, the $\Z_p$-modules $\Het^i(\Xbar, \Z_p(j))$ are finitely generated (\cite{milne}, Lemma V.1.11), and therefore their torsion subgroup is finite. For fixed $(i,j)$ this torsion subgroup is trivial for all but finitely many $p$ by a result of Gabber \cite{gabber}; as also remarked in (\cite{ctss}, p. 782, Remarque 1), in characteristic 0 the statement follows from comparison with the complex case.
\end{fact}

We also need a lemma on abelian groups.

\begin{lem}\label{konig}
Let $A$ be a finitely generated $\Z_p$-module equipped with a $\Z_p$-linear action of a group $G$. The group $(A\otimes\Q_p/\Z_p)^G$ is infinite if and only if $A^G$ contains an element of infinite order.
\end{lem}

\begin{dem} First of all, the finite torsion subgroup $T\subset A$ is a direct summand in $A$ which is also $G$-equivariant. Therefore after replacing $A$ by $A/T$ we may assume that $A$ is a free $\Z_p$-module. Then $A\otimes\Q_p/\Z_p$ is the $G$-equivariant direct limit of the finite groups $A/p^nA$ for all $n$, and each map $A/p^nA\to A\otimes\Q_p/\Z_p$ is injective with image equal to the $p^n$-torsion part of $A\otimes\Q_p/\Z_p$.

Now if $A^G$ contains an element of infinite order, the order of the groups $(A/p^nA)^G$ is unbounded as $n$ goes to infinity, and hence their direct limit $(A\otimes\Q_p/\Z_p)^G$ is infinite.

Conversely, assume that $(A\otimes\Q_p/\Z_p)^G$ is infinite.  Associate a graph to $(A\otimes\Q_p/\Z_p)^G$ whose vertices correspond to elements $a\in (A\otimes\Q_p/\Z_p)^G$ and two vertices $a$ and $a'$ are joined by an edge if $pa=a'$ or $pa'=a$. This graph is an infinite rooted tree in which each vertex has finite degree. Therefore by K\"onig's lemma  (see e.g. \cite{nath}, Theorem 3) it has an infinite path beginning at the root (which is the vertex corresponding to the element 0). This means that there exists an infinite sequence $(a_n)\subset (A\otimes\Q_p/\Z_p)^G$  such that $a_0=0$ and $pa_n=a_{n-1}$ for all $n>0$. In particular, $a_n$ has exact order $p^n$ in $(A\otimes\Q_p/\Z_p)^G$. As we may identify each $a_n$ with an element in $(A/p^nA)^G$, the sequence $(a_n)$  defines an element of infinite order in  $\limproj (A/p^nA)^G=A^G$.
\end{dem}

\begin{rema}\rm
As Wayne Raskind reminds us, it is possible to give a more traditional proof of the lemma, under the further assumption that $G$ is a compact topological group acting continuously on $A$ as above (which will be the case in our application). Indeed, since $(A\otimes_{\Z_p}\Q_p/\Z_p)^G$ is a $\Z_p$-module of finite cotype, its quotient modulo the maximal divisible subgroup $D$ is finite. Therefore we have to prove that $D$ is trivial if and only if $A^G$ is torsion. The latter condition is equivalent to the vanishing of the group $A^G\otimes_{\Z_p}\Q_p=(A\otimes_{\Z_p}\Q_p)^G$. By (\cite{tate}, Proposition 2.3), $D$ equals the image of the last map in the exact sequence of $\Z_p$-modules
$$
0\to A^G\to (A\otimes_{\Z_p}\Q_p)^G\to (A\otimes_{\Z_p}\Q_p/\Z_p)^G.
$$
Thus $D=0$ if and only if the map $A^G\to (A\otimes_{\Z_p}\Q_p)^G$ is surjective, which holds if and only if $(A\otimes_{\Z_p}\Q_p)^G=0$.
\end{rema}

\noindent{\em Proof of Proposition \ref{ABlemma}.\/} The finiteness of $\Het^i(\bar X, {\bf {\Q}}/{\Z}(j))^G$ is equivalent to the conjunction of the following two statements:\smallskip

\begin{enumerate}
\item[\rm (A')] $\Het^i(\bar X, {\Q}_p/\Z_p(j))^G$ is finite for all primes $p$.
\item[\rm (B')] $\Het^i(\bar X, {\Q}_p/{\Z}_p(j))^G=0$ for all but finitely many primes $p$.
\end{enumerate}\smallskip

We first prove (A) $\Leftrightarrow$ (A'). Consider the exact sequence of Galois modules
\begin{equation}\label{pex}
0 \to \Het^{i}(\bar X, {\Z}_p(j))/p^n \to \Het^i(\bar X, {\Z}/p^n{\Z}(j))  \to {}_{p^n}\Het^{i+1}(\bar X, {\Z}_p(j))\to 0
\end{equation}
coming from the long exact sequence associated with the multiplication-by $p^n$ map on $\Z_p(j)$; here ${}_{p^n}A$ denotes the $p^n$-torsion part of an abelian group $A$. By passing to the direct limit over all $n$, we obtain an exact sequence of $G$-modules
$$
0 \to \Het^{i}(\bar X, {\Z}_p(j))\otimes\Q_p/\Z_p \to \Het^i(\bar X, {\Q}_p/{\Z}_p(j))  \to \Het^{i+1}(\bar X, {\Z}_p(j))_{\rm tors}\to 0.
$$
As recalled above, the third group in this sequence is finite. By taking $G$-invariants we therefore see that (A') is equivalent to $(\Het^{i}(\bar X, {\Z}_p(j))\otimes\Q_p/\Z_p)^G$ being finite for all $p$. By Lemma \ref{konig} this is equivalent to $\Het^{i}(\bar X, {\Z}_p(j))^G$ being torsion for all $p$, which is the same as (A).

To finish the proof of the proposition, we prove (B) $\Rightarrow$ (B'). In view of Fact \ref{gabber} we may assume that the groups  $\Het^{i}(\bar X, {\Z}_p(j))$ and  $\Het^{i+1}(\bar X, {\Z}_p(j))$ are torsion free. Then exact sequence (\ref{pex}) yields isomorphisms
$$\Het^i(\bar X, {\Z}_p(j))/p^n \cong \Het^i(\bar X, {\Z}/p^n{\Z}(j))$$
for all $n$. Therefore we obtain exact sequences
\begin{equation}\label{pnex} 0 \to \Het^i(\bar X, {\Z}/p^{n-1}{\Z}(j)) \to \Het^i(\bar X, {\Z}/p^n{\Z}(j)) \to \Het^i(\bar X, {\Z}/p{\Z}(j)) \to 0\end{equation}
from tensoring the exact sequence
$$
0\to {\Z}/p^{n-1}{\Z}\to {\Z}/p^{n}{\Z}\to {\Z}/p{\Z}\to 0
$$
by the group $\Het^i(\bar X, {\Z}_p(j))$ which was also assumed to be torsion free. The sequences (\ref{pnex}) are $G$-equivariant, so after taking $G$-invariants a straightforward induction on $n$ shows that $H^i(\bar X, {\Z}/p^n{\Z}(j))^G$ vanishes for all $n$ provided it vanishes for $n=1$.
\enddem

\section{The oddity of cohomology}\label{oddity}

In this section we prove Theorem \ref{cohribet}, of which we take up the notation. We have to verify conditions (A) and (B) of Proposition \ref{ABlemma} in our situation. For condition (A) we prove a  vanishing result that generalizes (\cite{ribet}, Theorem 3).

\begin{prop}\label{vanish}
If $i$ is odd, we have $\Het^i(\Xbar, \Q_p(j))^G=0$ for all primes $p$.
\end{prop}

We adapt the proof of \cite{ribet} (itself based upon arguments of Imai \cite{imai} and Serre).
It uses the following fact from algebraic number theory:

\begin{lem}\label{riblem}
For every $p$ the largest subextension of $K|k$ unramified outside the primes dividing $p$ and infinity is obtained as the composite of $k(\mu_{p^\infty})$ with the largest subextension of $K|k$ unramified at all finite primes (which is a finite extension).
\end{lem}

\begin{dem}
This is the Lemma on p. 316 of \cite{ribet}.
\end{dem}

We also need the following consequence of the local monodromy theorem.

\begin{lem}\label{locmon} ${}$

\noindent $a)$ (Weak form) Fix $i,j$ and $p$.
There is a finite extension $k'|k$ such that every inertia subgroup in $\Gal(\ov k|k')$ associated with a prime not lying above $p$ acts unipotently on $\Het^i(\bar X, {\Q}_p(j))$.

\noindent $b)$ (Strong form) Fix $i$ and $j$.
There is a finite extension $k'|k$ such that for all primes $p$, every inertia subgroup in $\Gal(\ov k|k')$
associated with a prime not lying above $p$ acts unipotently on $\Het^i(\bar X, {\Q}_p(j))$.
\end{lem}

\begin{dem} Since $X$ extends to a smooth proper scheme over an open subscheme of the ring of integers of $k$, there are only finitely many conjugacy classes of inertia subgroups in $\Gal(\ov k|k)$ that act nontrivially on $\Het^i(\bar X, {\Q}_p(j))$, where $p$ is any prime number not dividing the residue characteristic associated with the inertia subgroup.

Fix first a prime $p$. By Grothendieck's local monodromy theorem (\cite{serretate}, Appendix) all inertia subgroups in $\Gal(\ov k|k)$ not associated with primes above $p$ act quasi-unipotently on $\Het^i(\bar X, \Q_p(j))$. Since up to conjugacy there are only finitely many that act nontrivially, their action becomes unipotent after replacing $k$ by a suitable finite extension. This yields $a)$.

The proof of $b)$ is similar, except that we use the strong version of the local monodromy theorem that yields an open subgroup of inertia acting unipotently on all $\Het^i(\bar X, \Q_p(j))$ for $p$ different from the residue characteristic. It is a consequence of de Jong's alteration theorem and the vanishing cycle spectral sequence (\cite{berthelot}, Proposition 6.3.2).
\end{dem}

\noindent{\em Proof of Proposition \ref{vanish}.\/}
The first step is, roughly,  to replace $K$ by $k(\mu_{p^\infty})$. This is achieved as follows. Assume $p$ is such that $\Het^i(\Xbar, \Q_p(j))^G \neq 0.$ By enlarging $k$ if necessary, we may assume $\mu_p\subset k$. The Galois group $\Gamma:=\Gal(\overline k|k)$ acts on $\Het^i(\Xbar, \Q_p(j))^G$ via its quotient $\Gamma/G$.   Choose a simple nonzero $\Gamma$-submodule $W$ of  $\Het^i(\Xbar, \Q_p(j))^G$. As $\Gamma/G$ is abelian and $W$ is simple, the elements of $\Gamma$ act semisimply on $W$. Therefore Lemma \ref{locmon} $a)$  implies that up to replacing $k$ by a finite extension we may assume that the action of $\Gamma$ on $W$ is unramified at all primes not dividing $p$. Then Lemma \ref{riblem} implies that, again up to replacing $k$ by a finite extension, the action of $\Gamma$ on $W$ factors through $\Gamma_p:=\Gal(k(\mu_{p^\infty})|k)$.

Let $D_v\subset \Gamma$ be a decomposition group of a prime $v$ of $k$ dividing $p$. Since $v$ is totally ramified in the extension $k(\mu_{p^\infty})|k$, the image of $D_v$ by the surjection $\Gamma\twoheadrightarrow \Gamma_p$ is the whole of $\Gamma_p$.  On the other hand, the abelian semisimple representation $W$ of $\Gamma$ restricts to a Hodge--Tate representation of $D_v$, by Hodge--Tate decomposition \cite{faltings} of the \'etale cohomology group $\Het^i(\Xbar, \Q_p(j))$. Therefore, by a theorem of Tate (\cite{serreab}, III, Appendix, Theorem 2), some open subgroup of $D_v$, and hence of $\Gamma$, acts on $W$ via the direct sum of integral powers of the $p$-adic cyclotomic character $\chi_p$. Replacing $k$ by a finite extension for the last time, we may assume that the whole of $\Gamma$ acts in this way. A Frobenius element $F_w$ at a prime $w$ of good reduction thus acts with eigenvalues that are integral powers of $\chi_p(F_w)=Nw$ (the cardinality of the residue field of $w$). But by the Weil conjectures as proven by Deligne, these eigenvalues should have absolute value $(Nw)^{i/2-j}$, a contradiction for odd $i$.
\enddem

In order to verify condition (B) of Proposition \ref{ABlemma}, we prove:

\begin{prop}\label{vanishmodp}
If $i$ is odd, we have $\Het^i(\Xbar, \Z/p\Z(j))^G=0$ for all but finitely many primes $p$.
\end{prop}

 \begin{dem} As in the proof of Proposition \ref{vanish}, we are allowed to replace $k$ by a finite extension throughout the proof. First we replace $k$ by the finite extension $k'|k$ given by Lemma \ref{locmon} $b)$ (this involves changing $K$ as well). Next, we replace $k$ by its maximal extension contained in $K$ in which no finite prime ramifies.

For all but finitely many $p$ the following conditions are all satisfied:
\begin{enumerate}
\item $p$ is unramified in $k$.
\item $\mu_p\not\subset k$.
\item $X$ has good reduction at the primes dividing $p$.
\item\label{free} The groups  $\Het^{i}(\bar X, {\Z}_p(j))$ and  $\Het^{i+1}(\bar X, {\Z}_p(j))$ are torsion free (see Fact \ref{gabber}).

\end{enumerate}

Assume now that there exist infinitely may primes $p$ satisfying the conditions above for which $\Het^i(\Xbar, \Z/p\Z(j))^G\neq 0$. We shall derive a contradiction.

For each such $p$ the nontrivial group $\Het^i(\Xbar, \Z/p\Z(j))^G$ carries an action of $\Gamma$ via the quotient $\Gamma/G$. The first step is to show that the restriction of this action to a simple $\Gamma$-submodule $W_p$ of $\Het^i(\Xbar, \Z/p\Z(j))^G$ factors through $\Gal(k(\mu_p)|k)$. Here we are following Ribet's argument in the proof of (\cite{ribet}, Theorem 2) closely. The image of the map $\Gamma\to \End_{\F_p}(W_p)$ giving the action of $\Gamma$ lies in $\End_\Gamma(W_p)$ as the action of $\Gamma$ on $W_p$ factors through its abelian quotient $\Gamma/G$. By Schur's lemma $\End_\Gamma(W_p)$ is a finite-dimensional division algebra over the finite field $\F_p$, hence a finite field $\F$ by Wedderburn's theorem. Thus the action of $\Gamma$ on $W_p$ is given by a character $\Gamma\to\F^\times$; in particular, it is semisimple. Since we have extended $k$ so that the conclusion of Lemma \ref{locmon} $b)$ holds, we know that every inertia group $I\subset\Gamma$ associated with a prime not dividing $p$ acts unipotently on $\Het^{i}(\bar X, {\Z}_p(j))$. Furthermore,  condition (\ref{free}) implies that we have an isomorphism $\Het^i(\Xbar, \Z/p\Z(j))\cong \Het^i(\Xbar, \Z_p(j))\otimes\Z/p\Z$ as in the proof of Proposition \ref{ABlemma}. Therefore $I$ acts on $\Het^i(\Xbar, \Z/p\Z(j))$ with eigenvalues congruent to 1 modulo $p$. In particular this applies to its action on $W_p$ which is semisimple, and therefore it must be trivial. Hence the action of $\Gamma/G$ on $W_p$ factors through its largest prime-to-$p$ quotient unramified outside infinity and the primes dividing $p$. By
Lemma \ref{riblem}, our assumption on $k$ implies that this quotient is $\Gal(k(\mu_p)|k)$. As we assumed $\mu_p\not\subset k$,  this is a nontrivial group isomorphic to $\F_p^\times$. Since $W_p$ is simple for the action of $\Gamma$, it must be 1-dimensional over $\F_p$, with $\Gamma$ acting by a power $\bar\chi_p^{n(p)}$ of the mod $p$ cyclotomic character $\bar\chi_p$.

Now by Serre's tame inertia conjecture (proven in \cite{caruso}, Theorem 5.3.1 -- in fact we only need the good reduction case which already follows from the work done in \cite{fola} and \cite{fome}), there exists a bound $N$ {\em independent of $p$} such that the integer $n(p)$ appearing in the above action satisfies $n(p)\leq N$. [In fact, Serre's conjecture concerns powers not of $\bar\chi_p$, but rather of the fundamental character denoted by $\theta_{p-1}$ in (\cite{serrepg}, \S 1). This character is associated with the Galois action on a $(p-1)$-st root of a uniformizer of a valuation of $k$  dividing $p$. However, by {\em loc. cit.}, Corollary to Proposition 8, we have $\bar\chi_p=\theta_{p-1}$ in case $p$ is  unramified in $k$, which we assumed.]

Now choose a place $w$ of $k$ not dividing $p$ where $X$ has good reduction, and let $Nw$ be the cardinality of its residue field. Under the isomorphism $\Gal(k(\mu_p)|k)\stackrel\sim\to \F_p^\times$ induced by $\bar\chi_p$ the Frobenius element $\bar F_w$ of $w$ in $\Gal(k(\mu_p)|k)$ corresponds to the class of $Nw$ modulo $p$. Therefore its action on $W_p$ is given by multiplication by $(Nw)^{n(p)}$, where $n(p)$ is the integer of the previous paragraph.
Now lift $\bar F_w$ to a Frobenius element $F_w$ of $w$ in $\Gamma$, and let $Q$ be the characteristic polynomial of its action on $\Het^i(\bar X,\Q_p(j))$; it is known that $Q$ has integral coefficients.  Moreover, by the Cayley--Hamilton theorem the minimal polynomial of $F_w$ acting on $\Het^i(\Xbar, \Z/p\Z(j))\cong \Het^i(\Xbar, \Z_p(j))\otimes\Z/p\Z$  divides the reduction of $Q$ modulo $p$, since $\Het^i(\Xbar, \Z_p(j))$ is torsion free by our assumption (\ref{free}). By the first part of the proof, the restriction of this action to a nonzero simple $\Gamma$-submodule $W_p\subset \Het^i(\Xbar, \Z_p(j))^G$  factors through $\Gal(k(\mu_p)|k)$ and corresponds to multiplication by $(Nw)^{n(p)}$ for some integer $n(p)$, so we conclude $Q((Nw)^{n(p)})\equiv 0$ modulo $p$. By the previous paragraph, this congruence holds for infinitely many $p$ but with $n(p)$ varying between 0 and a fixed bound $N$. Hence for some integer ${0\leq n(p)\leq N}$ we must have $Q((Nw)^{n(p)})=0$. But by the Weil conjectures proven by Deligne, we must then have $(Nw)^{n(p)}=(Nw)^{i/2-j}$, which is impossible for odd $i$. This gives our desired contradiction, so we finally conclude that $\Het^i(\Xbar, \Z/p\Z(j))^G$ vanishes for all but finitely many $p$.
\end{dem}

\begin{remas}\rm ${}$\smallskip

\noindent 1. In his proof of (\cite{ribet}, Theorem 2) which corresponds to the case $i=j=1$ here, Ribet used the Oort--Tate classification of finite group schemes at the point where we invoked Serre's tame inertia conjecture. On the other hand, instead of our final weight argument (which nicely parallels the proof of Proposition \ref{vanish}) he exploited the finiteness of global torsion on abelian varieties over $k$. In the general case we do not have a corresponding motivic finiteness theorem at our disposal.\smallskip

\noindent 2. Although the two statements are different in nature, there are remarkable similarities between the above proof and that of Faltings \cite{faltingsmordell} for his theorem that the height of abelian varieties over number fields is bounded in an isogeny class. Compare especially with the rendition by Deligne (\cite{delignefaltings}, Theorem 2.4).
\end{remas}

As already remarked, Theorem \ref{cohribet} is an immediate consequence of Propositions \ref{vanish} and \ref{vanishmodp} in view of Proposition \ref{ABlemma}.

\begin{rema}\label{remeven}\rm Theorem \ref{cohribet} does not hold for even degree cohomology. In fact, for projective space $\P^n_k$ over $k$ we have an isomorphism of Galois modules $\Het^{2i}(\P^n_{\kbar}, \Q/\Z(j)) \cong \Q/\Z(j-i)$ for all $0\leq i\leq n$ (see e.g. \cite{milne}, VI.5.6). As our $K$ contains all roots of unity by assumption, $G$ acts trivially on all $\Q/\Z(j-i)$ and therefore its invariants are infinite. However, the odd degree cohomology of projective space is trivial.

As Jean-Louis Colliot-Th\'el\`ene points out, the groups $\Het^{2i}(\Xbar, \Q/\Z(j))^G$ are infinite for any smooth projective $k$-variety $X$, essentially for the same reason as in the case $X=\P^n_k$. To see this, use Bertini's theorem to find linear subspaces in the ambient projective space intersecting $X_K$ in smooth connected closed $K$-subvarieties $H$ and $Y$, where $H$ has codimension $i\/$ in $X_K$ and $Y$ has dimension $i$. We may assume they are in general position. Their scheme-theoretic intersection is then a zero-dimensional reduced closed subscheme of length equal to the degree $d$ of $X$. Consider the cyclic subgroup of the Chow group $CH^i(X_K)$ generated by the class of $H$. Its image by the composite map
$$
\alpha:\, CH^i(X_K)\to CH^i(Y)\to  \Het^{2i}(\overline Y, \Q/\Z(j))^G\stackrel\sim\to \Q/\Z
$$
is $d(\Q/\Z)=\Q/\Z$, because the cycle map $CH^i(Y)\to  \Het^{2i}(\overline Y, \Q/\Z(j))$ factors through the degree map on $CH^i(Y)=CH_0(Y)$. But the composite map $\alpha$ factors through the group $\Het^{2i}(\Xbar, \Q/\Z(j))^G$ by functoriality of the cycle map, and therefore $\Het^{2i}(\Xbar, \Q/\Z(j))^G$ has an infinite quotient.

It would be interesting to know whether the module of $G$-invariants can still be infinite if we replace $\Het^{2i}(\Xbar, \Q/\Z(j))$ by its quotient modulo the image of the codimension $i$ cycle map.
\end{rema}

\section{Large variation and finiteness of cohomology in positive characteristic}\label{largevariation}

This section is devoted to the large variation condition of Definition \ref{large}.\smallskip

We begin with the {\em proof of Proposition \ref{cohposribet}.} Consider a base field $k=\F(C)$ of positive characteristic and a proper smooth $k$-variety $X$. Our task is to show the finiteness of the group $\Het^i(\Xbar, (\Q/\Z)'(j))^G$, where $G=\Gal(\kbar|k\overline\F)$. For this we are allowed to take finite extensions of $\F$, and therefore to find  a proper flat model $\CX$ of $X$  over $C$ for which there are rational points $c_1, c_2\in C$ satisfying the property in Definition \ref{large}. Also, we may assume $j=0$ as $G$ fixes all roots of unity in $\kbar$.

By Proposition \ref{ABlemma} and Remark \ref{remABlemma} it will suffice to prove:
\begin{itemize}
\item[\rm (A)] $\Het^i(\bar X, {\Q}_\ell)^G=0$ for all primes $\ell\not=p$.
\item[\rm (B)] $\Het^i(\bar X, {\Z}/\ell\Z)^G=0$ for all but finitely many primes $\ell\neq p$.
\end{itemize}

To show (A), notice that by the exact sequence
$$
1\to G\to\Gal(\bar k|k)\to \Gal(\overline\F|\F)\to 1
$$
the action of $\Gal(\bar k|k)$ on $\Het^i(\Xbar, \Q_\ell)^G$ factors
through $\Gal(\bar\F|\F)$, and hence the restrictions of the endomorphisms $\Frob_{c_1}, \Frob_{c_2}\in\End({\Het^i(\bar X, {\Q}_\ell)})$ to ${\Het^i(\bar X, {\Q}_\ell)^G}$ coincide by $\F$-rational\-ity of the points $c_1$ and $c_2$. But by the large variation assumption the $\Frob_{c_r}$ ($r=1,2$) have coprime characteristic polynomials, so this is only possible if $\Het^i(\bar X, {\Q}_\ell)^G=0$.

To prove (B), notice first that the elements $\Frob_{c_r}$ already act on the groups $\Het^i(\bar X, {\Z}_\ell)$ and the latter groups are torsion free for $\ell$ large enough (see Fact \ref{gabber}). In this case we may speak of the characteristic polynomials of the $\Frob_{c_r}$ on $\Het^i(\bar X, {\Z}_\ell)$. They are the same as on $\Het^i(\bar X, {\Q}_\ell)$ and are independent of $\ell$. As by assumption the characteristic polynomials of the $\Frob_{c_r}$ on $\Het^i(\bar X, {\Q}_\ell)$ are coprime, the same holds for their characteristic polynomials on  $\Het^i(\bar X, {\Z}/\ell\Z)$ for $\ell$ large enough. Indeed, the eigenvalues of the $\Frob_{c_r}$ on $\Het^i(\bar X, {\Z}_\ell)$ are algebraic integers (\cite{sga7}, Expos\'e XII, Th\'eor\`eme 5.2.2), so different eigenvalues can coincide modulo $\ell$ for only finitely many $\ell$. The end of the proof  of statement (B) is then the same as that of (A).
\enddem

\begin{remas}\label{largevarremas}\rm ${}$\smallskip

\noindent 1. The above argument, though much more elementary, is quite similar to the proof of Theorem \ref{cohribet}: we are comparing eigenvalues of two different Frobenius elements. The different nature of $\ell$-adic and $p$-adic weights guarantees `large variation' in the arithmetic setting.\smallskip

\noindent 2. An argument similar to that in Remark \ref{remeven} shows that large variation can only occur for odd degree cohomology groups.

\end{remas}

\begin{remas}\label{abvarremas}\rm ${}$\smallskip

\noindent 1. In the case where $X$ is an abelian variety and $\CX$ its N\'eron model over $C$, the large variation assumption for $i=1$ and with respect to the model $\CX$ has the following geometric reformulation: there exist two closed points $c_1, c_2\in C$ whose associated geometric fibres are abelian varieties over $\overline \F$ having no common simple isogeny factor. Indeed, if such a common isogeny factor exists, then by the K\"unneth formula its \'etale cohomology is a direct summand in the cohomology of the whole fibres. After extending the base field the inclusion of this direct summand becomes compatible with Frobenius, and therefore the characteristic polynomials of Frobenius have a common factor. Conversely, if we know that the group $H^1(\Xbar, \Q_\ell(1))$ does not have large variation, then we can identify a common nonzero Galois-invariant subspace in $H^1(\CX_{\bar c_1}, \Q_\ell(1))$ and $H^1(\CX_{\bar c_2}, \Q_\ell(1))$ after finite extension of $\F$. By semisimplicity of Frobenius we may extend the identification of the common subspaces to a Galois-equivariant morphism $H^1(\CX_{\bar c_1}, \Q_\ell(1))\to H^1(\CX_{\bar c_2}, \Q_\ell(1))$. But by the Tate conjecture for abelian varieties over finite fields such a morphism comes from a $\Q_\ell$-linear combination of maps $\CX_{\bar c_1}\to \CX_{\bar c_2}$ of abelian varieties. Such a map can only be nonzero if source and target have a common simple isogeny factor.\smallskip

\noindent 2. If the $k|\F$-trace of $X$ is nontrivial, it is not hard to check that all geometric fibres must have a common simple isogeny factor. Douglas Ulmer has communicated to us an argument showing that the converse holds when $X$ is an elliptic curve. It would be nice to know whether the converse statement holds in arbitrary dimension. This question is thoroughly studied in work in progress by Cadoret and Tamagawa \cite{cata}. They show that the converse would follow from a conjecture of Zarhin \cite{zarhinussr}, and they also prove the analogous statement over finitely generated infinite fields using a Hilbert irreducibility argument.\smallskip

\noindent 3. In any case, under the large variation assumption we obtain a purely cohomological proof of the prime-to-$p$ torsion part of Theorem \ref{posribet} that does not use the Lang-N\'eron theorem.
\end{remas}

\section{Proofs of the results on cycles}

In this section we prove Theorems \ref{main1} and \ref{main2}. We shall mostly adapt arguments used in the study of codimension 2 cycles on varieties over number fields, for which our main reference is \cite{ctcime}. There are two notable differences: our fields $K$ are of cohomological dimension 1, not 2 (which simplifies matters), whereas in characteristic 0 the \'etale cohomology groups with finite coefficients of the ring of integers of $K$ are not all finite (which complicates matters).\medskip

We begin with the {\em proof of Theorem \ref{main1}.} Denote the torsion subgroup of an abelian group $A$ by $A_{\rm tors}$. Over the algebraic closure $\kbar$ we have Bloch's Abel--Jacobi map
$$
CH^i(\Xbar)_{\rm tors}\to H^{2i-1}(\Xbar, \Q/\Z(i))
$$
which is injective for $i=2$ (\cite{ctcime}, th\'eor\`eme 4.3). It is moreover functorial in $\Xbar$, hence Galois-equivariant. So for $G=\Gal(\kbar|K)$ we have an injection
$$
CH^2(\Xbar)_{\rm tors}^G\hookrightarrow H^{3}(\Xbar, \Q/\Z(2))^G
$$
where the group on the right hand side is finite by Theorem \ref{cohribet}. It therefore suffices to show that the group $\ker(CH^2(X_K)\to CH^2(\Xbar))$ (which is torsion by the standard restriction-corestriction argument) has finite exponent and in fact trivial under the additional assumption on $\Het^3(\Xbar,\Z_\ell)$. This we do by a slight improvement of arguments of Colliot-Th\'el\`ene and Raskind \cite{ctra1}.

By (\cite{ctra1}, Proposition 3.6) we have an exact sequence
\begin{equation}\label{blochex}
H^1(K, K_2(\kbar(\Xbar))/H^0_{\rm Zar}(\Xbar, {\mathcal K}_2))\to \ker(CH^2(X_K)\to CH^2(\Xbar))\to H^1(K, H^1_{\rm Zar}(\Xbar, {\mathcal K}_2)).
\end{equation}
This exact sequence is obtained, following Bloch, by analyzing the commutative diagram of Gersten complexes
$$
\begin{CD}
 K_2(K(X_K)) @>>> \bigoplus_{x\in X_K^1} K(x)^\times @>>> \bigoplus_{x\in X_K^2} \Z \\
 @VVV @VVV @VVV \\
 K_2(\kbar(\Xbar))^G @>>> \left(\bigoplus_{x\in \Xbar^1} \kbar(x)^\times\right)^G @>>> \left(\bigoplus_{x\in \Xbar^2} \Z\right)^G
\end{CD}
$$
and using the fundamental fact proven by Quillen that the Gersten complexes compute the Zariski cohomology of the sheaf ${\mathcal K}_2$ (with $H^2_{\rm Zar}(X, {\mathcal K}_2)\cong CH^2(X)$ according to Bloch).

We now show that $H^1(K, K_2(\kbar(\Xbar))/H^0_{\rm Zar}(\Xbar, {\mathcal K}_2))=0$. A key point is that our field $K$ has cohomological dimension 1. Indeed, $K$ is the maximal cyclotomic extension of a number field, hence satisfies the assumptions of  (\cite{cogal}, II.3.3, Proposition 9). Now consider the exact sequence of $G$-modules
\begin{equation}\label{h0k2}
0\to H^0_{\rm Zar}(\Xbar, {\mathcal K}_2)\to K_2(\kbar(\Xbar))\to K_2(\kbar(\Xbar))/H^0_{\rm Zar}(\Xbar, {\mathcal K}_2)\to 0.
\end{equation}
By (\cite{ctra1}, Theorem 1.8) there is an exact sequence of Galois modules
$$
0\to T\to H^0_{\rm Zar}(\Xbar, {\mathcal K}_2)\to S\to 0,
$$
where $T$ is divisible as an abelian group and $S$ is a finite group which is the direct sum of the torsion subgroups in $H^2(\Xbar, \Z_\ell(2))$ for all $\ell$. In particular, $H^2(K, S)=0$ since $K$ has cohomological dimension 1. For a similar reason we have $H^2(K,T)=0$, as $T$ sits in an exact sequence of Galois modules
$$
0\to T_{\rm tors} \to T\to Q\to 0
$$
with  $Q$ uniquely divisible. Therefore we obtain $H^2(K, H^0_{\rm Zar}(\Xbar, {\mathcal K}_2))=0$.  On the other hand, ${\rm cd}(K)\leq 1$ also implies $H^1(K,K_2(\kbar(\Xbar)))=0$ in view of  (\cite{CT}, Corollary 1, p. 11). Thus the vanishing of $H^1(K, K_2(\kbar(\Xbar))/H^0(\Xbar, {\mathcal K}_2))$ follows from the long exact cohomology sequence of (\ref{h0k2}).

Finally, the group $H^1(K,H^1_{\rm Zar}(\Xbar, {\mathcal K}_2))$ is of finite exponent when  $H^2(X, \calo_X)=0$ and trivial under the additional assumption that $\Het^3(\Xbar, \Z_\ell)$ is torsion free for all $\ell$ (\cite{ctra1}, Proposition 3.9 b) and d)). The theorem thus results from exact sequence (\ref{blochex}).

\begin{rema}\rm
Though the group $H^1(K,H^1_{\rm Zar}(\Xbar, {\mathcal K}_2))$ has finite exponent as recalled above, it is in general infinite over our field $K$. However, this does not contradict the conjectured finiteness of $CH^2(X_K)_{\rm tors}$. Indeed, by analyzing the diagram of Gersten complexes in the above proof further, one may continue exact sequence (\ref{blochex}) as
$$
 \ker(CH^2(X_K)\to CH^2(\Xbar))\to H^1(K, H^1_{\rm Zar}(\Xbar, {\mathcal K}_2))\to H^2(K, K_2(\kbar(\Xbar))/H^0_{\rm Zar}(\Xbar, {\mathcal K}_2))
$$
(see \cite{ctra1}, Proposition 3.6). By similar arguments as above, we have
$$
H^2(K, K_2(\kbar(\Xbar))/H^0_{\rm Zar}(\Xbar, {\mathcal K}_2))\cong H^2(K, K_2(\kbar(\Xbar))).
$$
Here the group on the right hand side is usually large. In fact, by similar reduction arguments as in the proof of (\cite{CT}, Theorem B) one reduces its study to that of the group $H^2({\rm Gal}(L|K), K_2(L(X_L)))$ for a finite cyclic extension $L|K$. By periodicity of the cohomology of cyclic groups, this group is related to the cokernel of the norm map $K_2(L(X_L))\to K_2(K(X))$ which is large for a field of cohomological dimension $>2$ such as $K(X)$. On the other hand, it seems difficult to analyze the kernel of the map $H^1(K, H^1_{\rm Zar}(\Xbar, {\mathcal K}_2))\to H^2(K, K_2(\kbar(\Xbar))$.\end{rema}

We now begin the proof of Theorem \ref{main2}. The first point is:

\begin{prop}\label{main2.1}
Under the assumptions of Theorem \ref{main2} the prime-to-$p$ torsion subgroup of $CH^2(X_K)$ has finite exponent.
\end{prop}

The proof hinges on the following lemma.

\begin{lem}
Let $Y$ be a smooth projective variety over an algebraically closed field $F$ of characteristic $p>0$ that is liftable to characteristic 0. If moreover $H^2_{\rm Zar}(Y,\calo_Y)=0$,  the cokernel of the Chern class map
$$
\Pic(Y)\otimes\Q_\ell/\Z_\ell\stackrel{c_Y}\longrightarrow \Het^2(Y, \Q_\ell/\Z_\ell(1))
$$
is finite for all $\ell\neq p$, and zero for all but finitely many $\ell$ .

\end{lem}

Another way of phrasing the conclusion of the lemma is that the prime-to-$p$ torsion part of the Brauer group ${\rm Br}(Y)$ is finite (cp. \cite{gb2}, Theorem 3.1).\medskip

\begin{dem}
By the liftability assumption $Y$ arises as the special fibre of a smooth projective flat morphism $\CY\to S$, where $S$ is the spectrum of a complete discrete valuation ring of characteristic 0 and residue field $F$. If $\bar\eta$ denotes a geometric point over the generic point $\eta$ of $S$, we have a commutative diagram
\begin{equation}\label{diag}
\begin{CD}
\Pic(\CY_{\bar\eta}) @<<< \Pic(\CY) @>>> \Pic(Y) \\
@VV{c_{\CY_{\bar\eta}}}V @VV{c_{\CY}}V  @VV{c_Y}V\\
\Het^2(\CY_{\bar\eta}, \Q_\ell/\Z_\ell(1)) @<\cong<< \Het^2(\CY, \Q_\ell/\Z_\ell(1)) @>\cong>> \Het^2(Y, \Q_\ell/\Z_\ell(1))
\end{CD}
\end{equation}
where the horizontal maps are pullbacks and the lower horizontal maps are isomorphisms by the proper and smooth base change theorems.

By semi-continuity of coherent cohomology (\cite{hartshorne}, Chapter III, Theorem 12.8), the assumption
$H^2_{\rm Zar}(Y, \calo_{Y})=0$ forces  $H^2_{\rm Zar}(\CY_{\eta}, \calo_{\CY_{\eta}})=0$, and hence also $H^2_{\rm Zar}(\CY_{\bar\eta}, \calo_{\CY_{\bar\eta}})=0$ as the coherent cohomology of varieties behaves well with respect to extension of the base field. Since $\kappa(\bar\eta)$ has characteristic 0, the vanishing of $H^2_{\rm Zar}(\CY_{\bar\eta}, \calo_{\CY_{\bar\eta}})$ implies that the cokernel of the map $$\Pic(\CY_{\bar\eta})\otimes \Q_\ell/\Z_\ell\to \Het^2(\CY_{\bar\eta}, \Q_\ell/\Z_\ell(1))$$  is finite for all $\ell\neq p$ and zero for almost all $\ell$ (see e.g. \cite{ctra1}, Proposition 2.11). Thus if we show that, after identifying the groups $\Het^2(\CY_{\bar\eta}, \Q_\ell/\Z_\ell(1))$ and $\Het^2(Y, \Q_\ell/\Z_\ell(1))$ by means of the base change isomorphisms in the diagram, we have $\im(c_Y)\supset\im(c_{\CY_{\bar\eta}})$, the proposition will follow.

To verify this claim,  pick a class $\bar\alpha\in \Pic(\CY_{\bar\eta})$. There is a finite extension $L|\kappa(\eta)$ and a class $\alpha_L\in \Pic(\CY_L)$ mapping to $\bar\alpha$ under the pullback map $\Pic(\CY_L)\to \Pic(\CY_{\bar\eta})$, where $\CY_L:=\CY\times_S\Spec(L)$. Denote by $S'$ the normalization of $S$ in $L$ and by $\CY'$ the base change $\CY\times_SS'$. Since the map $S'\to S$ comes from a totally ramified extension of complete discrete valuation rings, the special fibre of $\CY'$ is still isomorphic to $Y$. As $\CY_L$ identifies with an open subscheme of $\CY'$, the pullback map $\Pic(\CY')\to\Pic(\CY_L)$ is surjective, and hence $\alpha_L$ comes from a class $\alpha'\in\Pic(\CY')$. This class in turn pulls back to $\alpha_Y\in \Pic(Y)$ on the special fibre. By construction, $c_Y(\alpha_Y)$ corresponds to $c_{\CY_{\bar\eta}}(\bar\alpha)$ under the base change isomorphism $\Het^2(\CY_{\bar\eta}, \Q_\ell/\Z_\ell(1)) \cong \Het^2(Y, \Q_\ell/\Z_\ell(1))$ coming from $\CY_{\bar\eta}=\CY_{\bar\eta}'\to \CY'\leftarrow Y$. But then these classes also correspond under the base change isomorphism coming from $\CY_{\bar\eta}\to \CY\leftarrow Y$,  in view of the commutative diagram
$$
\begin{CD}
\Het^2(\CY_{\bar\eta}, \Q_\ell/\Z_\ell(1)) @<\cong<< \Het^2(\CY, \Q_\ell/\Z_\ell(1)) @>\cong>> \Het^2(Y, \Q_\ell/\Z_\ell(1)) \\
@VV{\id}V @VVV  @VV{\id}V\\
\Het^2(\CY_{\bar\eta}, \Q_\ell/\Z_\ell(1)) @<\cong<< \Het^2(\CY', \Q_\ell/\Z_\ell(1)) @>\cong>> \Het^2(Y, \Q_\ell/\Z_\ell(1)).
\end{CD}
$$

\end{dem}

\begin{rema}\rm  The inclusion  $\im(c_Y)\supset\im(c_{\CY_{\bar\eta}})$ in the above proof is in fact an equality. The reverse inclusion follows from the commutative diagram (\ref{diag}) above, together with the surjectivity of the pullback map  $\Pic(\CY)\to\Pic(Y)$. This surjectivity is one of Grothendieck's deformation results exposed in  (\cite{fga}, \S 6) or (\cite{illusie}, \S 8.5); it uses the assumption $H^2_{\rm Zar}(Y,\calo_Y)=0$.
\end{rema}

\noindent{\em Proof of Proposition \ref{main2.1}.} Given the lemma above, the proof is a minor modification of that of Theorem \ref{main1}. We review the main steps. As in the cited proof, but using the large variation assumption and Proposition \ref{cohposribet} instead of Theorem \ref{cohribet}, we reduce to proving that the prime-to-$p$ torsion in $\ker(CH^2(X_K)\to CH^2(\Xbar))$ has finite exponent. This we again do using exact sequence (\ref{blochex}). In exactly the same way as above, we obtain that the group $H^2(K, H^0_{\rm Zar}(\Xbar, {\mathcal K}_2))$ has no prime-to-$p$ torsion, the only difference in the argument being that when applying (\cite{ctra1}, Theorem 1.8) we only obtain prime-to-$p$ divisibility of the group $T$. Next, we also get that $H^1(K,K_2(\kbar(\Xbar)))$ has no prime-to-$p$ torsion by applying  (\cite{CT}, Corollary 1). (Note that the result is only stated there in characteristic 0 because the proof relies on a result of Bloch, cited as {\bf B6} in \cite{CT}, that was only known in characteristic 0 at the time. However, shortly afterwards Suslin extended Bloch's statement to arbitrary characteristic in \cite{suslin}, Theorem 3.5.) From these facts we conclude as above that $H^1(K, K_2(\kbar(\Xbar))/H^0(\Xbar, {\mathcal K}_2))$ has no prime-to-$p$ torsion.

It remains to show that the prime-to-$p$ torsion in $H^1(K,H^1_{\rm Zar}(\Xbar, {\mathcal K}_2))$ is of finite exponent under our assumptions. For this we have to extend (\cite{ctra1}, Proposition 3.9 b)) to positive characteristic. The cited result is an immediate consequence of (\cite{ctra1}, Theorem 2.12) which is again only stated in characteristic 0. However, by (\cite{ctra1}, Remark 2.14) the prime-to-$p$ version of the said theorem holds in positive characteristic if one knows that the prime-to-$p$ torsion part of the Brauer group ${\rm Br}(Y)$ is finite. This is the content of the lemma above.
\enddem

\begin{prop}\label{main2.2}
 Let $K=\overline \F(C)$ be as in Theorem \ref{posribet}, and let $\ell$ be a prime different from the characteristic $p$. For every smooth proper $K$-variety $X_K$ satisfying $H^2_{\rm Zar}(X_K,\calo_{X_K})=0$ the $\ell$-primary torsion subgroup $CH^2(X_K)\{\ell\}$ is of finite cotype.
\end{prop}

Recall that an $\ell$-primary torsion abelian group $A$ is of finite cotype if for all $i>0$ the multiplication-by-$\ell^i$ map on $A$ has finite kernel. This is equivalent to the $\Q_\ell/\Z_\ell$-dual being a finitely generated $\Z_\ell$-module.\medskip

\begin{dem}
We adapt an argument that is at the very end of both \cite{ctcime} and \cite{salb}. We may find a smooth affine $\overline\F$-curve $\overline U\subset C_{\overline \F}$  so that $X_K\to\Spec K$ extends to a smooth proper flat morphism $\CX\to \overline U$.

We first prove that the restriction map on torsion subgroups $CH^2(\CX)_{\rm tors}\to CH^2(X_K)_{\tors}$ is surjective. Localization in Gersten theory gives an exact sequence
\begin{equation}\label{gloc}
H^1_{\rm Zar}(X_K, {\mathcal K}_2)\stackrel\partial\to\bigoplus_{P\in \overline U_0}\Pic (\CX_P)\to CH^2(\CX)\to CH^2(X_K)\to 0.
\end{equation}
Since $\Q$ is flat over $\Z$, the sequence remains exact after tensoring with $\Q$, so it will suffice to prove that the map $CH^2(\CX)\otimes\Q\to CH^2(X_K)\otimes\Q$ is an isomorphism, or equivalently that the map $\partial\otimes\Q$ is surjective. By a norm argument we may check the latter statement after passing to a finite extension $L|K$ (which entails the replacement of $\overline U$ by a finite cover). Thus we may assume that $X_K(K)\neq\emptyset$ and the composite map $\Pic(X_K)\to \Pic(X_{\overline K})\to {\rm NS}(X_{\overline K})$ is surjective,  where NS denotes, as usual, the N\'eron--Severi group.  Furthermore, by shrinking $\overline U$ if necessary we may assume using semi-continuity of coherent cohomology (\cite{hartshorne}, Chapter III, Theorem 12.8) that $H^2_{\rm Zar}(\CX_{P}, \calo_{\CX_{P}})=0$ for each closed point $P\in \overline U_0$.  In this situation, and under the further assumption that $\Pic(\overline U)=0$,  the composite map
\begin{equation}\label{CTra}
H^1_{\rm Zar}(X_K, {\mathcal K}_2)\stackrel\partial\to\bigoplus_{P\in \overline U_0}\Pic (\CX_{P})\to\bigoplus_{P\in \overline U_0}{\rm NS} (\CX_{P})
\end{equation}
is surjective according to Lemma 3.2 of \cite{ctra2} (there is also a characteristic zero assumption in the cited lemma but it is not used). In our case the assumption $\Pic(\overline U)=0$ is not satisfied, but at least $\Pic(\overline U)$ is torsion since $\overline U$ is a smooth affine curve over $\overline\F$. Under this weaker assumption the proof of (\cite{ctra2}, Lemma 3.2) yields that the composite map (\ref{CTra}) has torsion cokernel (only the last three lines have to be modified in a straightforward way). Thus the composite map
$$
H^1_{\rm Zar}(X_K, {\mathcal K}_2)\otimes\Q\stackrel{\partial\otimes\Q}\longrightarrow\bigoplus_{P\in \overline U_0}\Pic (\CX_{P})\otimes\Q\to\bigoplus_{P\in \overline U_0}{\rm NS} (\CX_{P})\otimes\Q
$$
is surjective. But here the second map is an isomorphism because the groups $\Pic^0(\CX_{P})$ are torsion as well. This shows the surjectivity of the map $\partial\otimes\Q$, and hence the surjectivity of the map $CH^2(\CX)_{\rm tors}\to CH^2(X_K)_{\tors}$. In particular, we have a surjective map $CH^2(\CX)\{\ell\}\to CH^2(X_K)\{\ell\}$ on $\ell$-primary torsion.

Therefore to finish the proof it will suffice to prove that the group $CH^2(\CX)\{\ell\}$ is of finite cotype. By Bloch--Ogus theory and the Merkurjev--Suslin theorem we have a surjection
$$
H^1_{\rm Zar}(\CX, \R^2\pi_*\Q_\ell/\Z_\ell(2))\twoheadrightarrow CH^2(\CX)\{\ell\}
$$
as well as an injection
$$
H^1_{\rm Zar}(\CX, \R^2\pi_*\Q_\ell/\Z_\ell(2))\hookrightarrow\Het^3(\CX, \Q_\ell/\Z_\ell(2))
$$
where $\pi:\, \CX_{\mbox{\scriptsize\'et}}\to\CX_{\rm Zar}$ is the change-of-sites map (see e.g. \cite{ctcime}, \S 3.2). Here $\CX$ is a smooth variety over the algebraically closed field $\overline \F$, and therefore the groups  $\Het^i(\CX, \Q_\ell/\Z_\ell(j))$ are of finite cotype for all $i,j$.
\end{dem}

\noindent{\em Proof of Theorem \ref{main2}.\/} We have to prove that for all $\ell\neq p$ the $\ell$-primary torsion subgroups $CH^2(X_K)\{\ell\}$ are finite and that they are trivial for all but finitely many $\ell$. The latter statement follows from the fact that the whole prime-to-$p$ torsion subgroup in $CH^2(X_K)$ has finite exponent (Proposition \ref{main2.1}). Now Propositions \ref{main2.1} and \ref{main2.2} together imply that for each fixed $\ell\neq p$ the group $CH^2(X_K)\{\ell\}$ is finite.
\enddem

\begin{rema}\rm The above proofs show that the conclusion of Theorem \ref{main2} also holds for the prime-to-$p$ torsion in $CH^2(X)$ instead of $CH^2(X_K)$, even without assuming large variation. One only has to replace the application of Proposition \ref{cohposribet} by a weight argument.
\end{rema}

\end{document}